\documentclass[11pt,reqno]{amsart}
\usepackage{mathptmx}
\usepackage[latin9]{inputenc}
\usepackage{color}
\usepackage{amsthm}
\usepackage{amsbsy}
\usepackage{amssymb}
\usepackage[numbers]{natbib}

\makeatletter

\DeclareRobustCommand{\greektext}{%
  \fontencoding{LGR}\selectfont\def\encodingdefault{LGR}}
\DeclareRobustCommand{\textgreek}[1]{\leavevmode{%
  \IfFileExists{grtm10.tfm}{}{\fontfamily{cmr}}\greektext #1}}
\DeclareFontEncoding{LGR}{}{}
\DeclareTextSymbol{\~}{LGR}{126}

\numberwithin{equation}{section}
\numberwithin{figure}{section}

\usepackage{amsfonts}\usepackage{amsthm}\usepackage{enumerate}\usepackage{color}
\usepackage{graphicx}
\usepackage{epstopdf}

\theoremstyle{definition}

\numberwithin{equation}{section}

\makeatother

\begin{document}
\setcounter{page}{1}

\vspace*{1cm}

\title[Problem Structures in Superiorization]{Problem Structures in the Theory and Practice of Superiorization}

\author{Gabor T. Herman \vspace*{-0.6cm}
}

\maketitle
\begin{center}
\textit{\footnotesize{}Computer Science Ph.D. Program, The Graduate
Center, City University of New York, New York, NY 10016, USA}
\par\end{center}{\footnotesize \par}

\vskip 4mm \textbf{\footnotesize{}Abstract.}{\footnotesize{} The
purpose of this short paper is to identify the mathematical essence
of the superiorization methodology. This methodology has been developed
in recent years while attempting to solve specific application-oriented
problems. Consequently, superiorization is often presented using the
terminology of such problems. A more general approach is provided
here by discussing ideas related to superiorization in terms of an
abstract mathematical concept, referred to as a problem structure.}{\footnotesize \par}

\noindent {\footnotesize{}\vskip 1mm }\textbf{\footnotesize{}Keywords.}{\footnotesize{}
Superiorization; Optimization; Feasibility-seeking; Proximity function;
Perturbation resilience; Algorithm.}{\footnotesize \par}

\noindent {\footnotesize{}\vskip 1mm }\textbf{\footnotesize{}2010
Mathematics Subject Classification.}{\footnotesize{} 15A29, 65F22,
65K10, 90C26.}{\footnotesize \par}

\global\long\def\thefootnote{}
 \footnotetext{Corresponding author: Gabor T. Herman. 

E-mail addresses: gabortherman@yahoo.com. 

Received xxxxxx xx, 2019; Accepted xxxxxx xx, 20xx. }

\dedicatory{In the spirit of Mel Brook's \textit{History of the World Part 1}
that contains the following conversation: ``Occupation?'' ``Stand-up
Philosopher`` ``What?'' ``A Stand-up Philosopher!; I coalesce
the vapor of human experience into a viable and logical comprehension.''
``Oh! A Bullshit Artist!?!!''}

\section{Introduction}

The term ``superiorization'' (in the sense as it is used here) first
appeared about a decade ago and there are now about 100 publications
on the topic \citep{Cen19}. In view of this recent and rapid development,
the terminology in the literature is far from settled; here we adopt
a terminology based on \citep{Her12,CGHH19}.

Another aspect of the superiorization methodology is that it has been
driven by problems in various applications. In many of the publications,
the mathematical essence of superiorization got intermixed with the
specifics of such an application. Our purpose here is to present superiorization
in an abstract mathematical manner that does not rely on concepts
from the application areas. To avoid potential confusion between the
abstract theory and its practical realizations, we do not provide
examples to illustrate the abstract definitions; many such examples
appeared in the literature, see \citep{Cen19} and, in particular,
\citep{Her12}.

To achieve our purpose, we make use of the concept of a ``problem
structure.'' As will be seen, problem structures are rather minimal
abstract structures that are nevertheless rich enough for a discussion
of the superiorization methodology. Problem structures are formally
and completely defined in the next section; no additional information
should be assumed because of the name ``problem structure.'' That
name has been selected based on the historical applications of superiorization
and has nothing to do with its mathematical essence, which is the
topic of the current paper.

A final introductory comment: Some of the more recent literature distinguishes
between weak and strong superiorization \citep{Cen15}, with ``strong
superiorization'' corresponding to what was called simply ``superiorization''
in earlier papers on the topic, such as \citep{Her12}. The meaning
of ``superiorization'' in the current presentation is also ``strong
superiorization.''

\section{Formal Definitions}

We use \textcolor{black}{$\mathbb{\mathbb{R}}$} to denote the set
of real numbers and $\mathbb{R}_{+}$ to denote the set of nonnegative
real numbers.

A \textit{\textcolor{black}{problem structure}} is a triple $\left\langle \Omega,\mathbb{T},\mathcal{P}r\right\rangle $,
where
\begin{itemize}
\item $\Omega$ and $\mathbb{T}$ are nonempty sets and
\item $\mathcal{P}r$ is a function\textcolor{black}{{} on $\mathbb{T}$ is
such that, for every $T\in\mathbb{T}$, $\mathcal{P}r_{T}:\Omega\rightarrow\mathbb{R}_{+}$
(so $\mathcal{P}r_{T}$ maps $\Omega$ into }$\mathbb{R}_{+}$\textcolor{black}{).}
\end{itemize}
A \textit{\textcolor{black}{targeted problem structure}} is a quadruple
$\left\langle \Omega,\mathbb{T},\mathcal{P}r,\phi\right\rangle $,
where $\left\langle \Omega,\mathbb{T},\mathcal{P}r\right\rangle $
is a problem structure and $\phi$ is a function from $\Omega$ into
\textcolor{black}{$\mathbb{R}$.}

\textcolor{black}{Now we give some intuition behind these formal definitions.
We think of $\mathbb{T}$ as a problem set for which the complete
specification of each problem $T\in\mathbb{T}$ is provided by the
function $\mathcal{P}r_{T}$ that assigns a nonnegative real value
to every element $\boldsymbol{x}$ of $\Omega$. That value indicates
how undesirable $\boldsymbol{x}$ is as a solution to $T$. }$\mathcal{P}r$
is often referred to as a \textit{proximity function}. The function
$\phi$ assigns a \textcolor{black}{real value to every element $\boldsymbol{x}$
of $\Omega$.} \textcolor{black}{That value indicates the prior undesirability
of $\boldsymbol{x}$; that is its undesirability without considering
any specific problem $T\in\mathbb{T}$. }$\phi$ is often referred
to as a \textit{target function} \citep{CGHH19}. We emphasize, again,
that these intuitive hints are provided only so that the reader gets
a glimpse of the motivation behind the formal definitions; as far
as the mathematics is concerned, the definitions are complete without
these hints.

One issue that is not relevant from the point of view of providing
mathematical definitions, but is essential in practical applications,
is the computability of the functions used in the definitions. Such
computability will be assumed in this paper without further comments.
(Thus, we assume the existence and availability of computer code that
for any \textcolor{black}{$T\in\mathbb{T}$ and $\boldsymbol{x}\in\Omega$,
calculates $\mathcal{P}r_{T}\left(\boldsymbol{x}\right)$.)}

For a (fixed) problem structure $\left\langle \Omega,\mathbb{T},\mathcal{P}r\right\rangle $,
a $T\in\mathbb{T}$, an $\varepsilon\in\mathbb{R}_{+}$ and a sequence
$R=\left(\boldsymbol{x}^{k}\right)_{k=0}^{\infty}$ of elements of
$\Omega$, we use $O\left(T,\varepsilon,R\right)$ to denote the element
$\boldsymbol{x}\in\Omega$ that has the following properties: $\mathcal{P}r_{T}(\boldsymbol{x})\leq\varepsilon$\emph{
}and there is a nonnegative integer $K$ such that $\boldsymbol{x}^{K}=\boldsymbol{x}$
and, for all nonnegative integers $k<K$, $\mathcal{P}r_{T}\left(\boldsymbol{x}^{k}\right)>\varepsilon$.
Clearly, if there is such an $\boldsymbol{x}$, then it is unique.
If there is no such $\boldsymbol{x}$, then we say that $O\left(T,\varepsilon,R\right)$
is \textit{undefined,} otherwise we say that it is \textit{defined}.
An \emph{algorithm }\textbf{$\mathbf{P}$ }for a problem structure
$\left\langle \Omega,\mathbb{T},\mathcal{P}r\right\rangle $ and a
set $\Delta$ such that $\Omega\subseteq\Delta$ assigns to each problem
$T\in\mathbb{T}$ an operator $\mathbf{P}_{T}:\Delta\rightarrow\Omega$.
For any \emph{initial point} $\boldsymbol{x}\in\Omega$, \textbf{$\mathbf{P}_{T}$}
produces the infinite sequence $\left(\left(\mathbf{P}_{T}\right)^{k}\boldsymbol{x}\right)_{k=0}^{\infty}$
of elements of $\Omega$. Intuitively, given an algorithm $\mathbf{P}$
for a problem structure $\left\langle \Omega,\mathbb{T},\mathcal{P}r\right\rangle $,
a problem $T\in\mathbb{T}$, an $\varepsilon\in\mathbb{R}_{+}$ and
an $\boldsymbol{x}\in\Omega$, we may think of $O\left(T,\varepsilon,\left(\left(\mathbf{P}_{T}\right)^{k}\boldsymbol{x}\right)_{k=0}^{\infty}\right)$
as the $\varepsilon$-\textit{output} for problem $T$ of algorithm
$\mathbf{P}$ when initialized at $\boldsymbol{x}$. Such an output
$\boldsymbol{y}$ may be undefined, but if it is defined, then $\boldsymbol{y}\in\Omega$
and $\mathcal{P}r_{T}(\boldsymbol{y})\leq\varepsilon$.

We are now in position to say something about the nature of superiorization.
Intuitively, what we wish to do is to design a methodology that for
an algorithm (defined as above) will produce a superiorized version
of it; meaning a version whose performance is as good as that of the
original algorithm from the point of view of the proximity function
but is better from the point of view of the target function. We now
make this very rough statement mathematically somewhat more precise.

Consider the infinite sequence $\left(\left(\mathbf{P}_{T}\right)^{k}\boldsymbol{x}\right)_{k=0}^{\infty}$
produced by \textbf{$\mathbf{P}_{T}$} for the initial point $\boldsymbol{x}\in\Omega$.
If we define, for all $k\geq0$, $\boldsymbol{x}^{k}=\left(\mathbf{P}_{T}\right)^{k}\boldsymbol{x}$,
then we have that, for all $k\geq0$, $\boldsymbol{x}^{k+1}=\mathbf{P}_{T}\boldsymbol{x}^{k}$.
In superiorization, this sequence $R=\left(\boldsymbol{x}^{k}\right)_{k=0}^{\infty}$
is perturbed to get a new sequence $S=\left(\boldsymbol{x}^{k}\right)_{k=0}^{\infty}$
as follows: $\boldsymbol{x}^{0}=\boldsymbol{x}$ and, for all $k\geq0$,
$\boldsymbol{x}^{k+1}=\mathbf{P}_{T}\mathbf{S}\boldsymbol{x}^{k}$,
where $\mathbf{S}:\Omega\rightarrow\Delta$. In this general discussion
we do not get into the details of defining the superiorizing operator
$\mathbf{S}$, but we note that, in any case, $\boldsymbol{x}^{k}\in\Omega$,
for all $k\geq0$ in the new sequence $S$. One rather-strong mathematical
restatement of the wish expressed in the previous paragraph is that
$\mathbf{S}$ ought to be chosen so that whenever $O\left(T,\varepsilon,R\right)$
is defined, then $O\left(T,\varepsilon,S\right)$ is also defined
and $\phi\left(O\left(T,\varepsilon,S\right)\right)\leq\phi\left(O\left(T,\varepsilon,R\right)\right)$.
(In usual practice we would like the target value for $S$ to be much
smaller than the target value for $R$, but less-than-or-equal-to
is good enough for our mathematical discussion.) 

We are now in position to provide, within our general framework, a
pseudocode for a skeleton of the superiorized version of an algorithm
$\mathbf{P}$. This skeleton is valid for a number of specific situations,
the lines in the following pseudocode appear in both (of the pseudocodes
of) Algorithms 1 and 2 in \citep{CGHH19}. \bigskip{}

\textbf{Skeleton of the superiorized version of algorithm }$\mathbf{P}$\\

\textbf{set} $k=0$ 

\textbf{set} $\boldsymbol{x}^{k}=x$ 

$\cdots$ 

\textbf{repeat} 

$\qquad$\textbf{set} $n=0$ 

$\qquad$\textbf{set} $\boldsymbol{x}^{k,n}=\boldsymbol{x}^{k}$ 

$\qquad$\textbf{while} $n<N$ 

$\qquad\qquad\;\,$$\cdots$

$\qquad\qquad\;\,$\textbf{$\cdots$}

$\qquad$\textbf{set} $\boldsymbol{x}^{k+1}=\mathbf{P}_{T}\boldsymbol{x}^{k,N}$ 

$\qquad$\textbf{set} $k=k+1$\\
\\
We now comment on this skeleton and discuss the way(s) it may be completed
by insertion of extra lines to provide a pseudocode for a complete
superiorized version of algorithm $\mathbf{P}$. 

The $\boldsymbol{x}$ in the second line of the skeleton is the initial
point for the algorithm. 

The line(s) to be inserted above \textbf{repeat} may initialize some
variables to be used inside the \textbf{while} loop. 

The skeleton indicates that the superiorization operator $\mathbf{S}$
is implemented in a user-specified number, $N$, stages. Each stage
is a small perturbation to get from an $\boldsymbol{x}^{k,n}$ in
$\Delta$ to an $\boldsymbol{x}^{k,n+1}$ in $\Delta$, with $\boldsymbol{x}^{k,0}=\boldsymbol{x}^{k}$,
and with the combined effect of producing $\boldsymbol{x}^{k,N}=\mathbf{S}$$\boldsymbol{x}^{k}$,
also in $\Delta$. From this we see that applying the step following
the \textbf{while} loop results in $\boldsymbol{x}^{k+1}=\mathbf{P}_{T}\mathbf{S}\boldsymbol{x}^{k}$
(which is in $\Omega)$, as promised prior to giving the pseudocode
for the skeleton. 

We next discuss the nature of the code to be inserted into the skeleton
in the \textbf{while} loop, namely after the line \textbf{while} $n<N$
and before the line \textbf{set} $\boldsymbol{x}^{k+1}=\mathbf{P}_{T}\boldsymbol{x}^{k,N}$.
The first time that code is entered, it is the case that $n=0$. In
order to keep a record of where we are in the the execution of the
$N$ stages, there should be a to-be-executed statement of the form
$\mathbf{set}$ $n=n+1$ in the inserted code. After the $N$th execution
of this statement we will have $n=N$ and this indicates the last
execution of the code in the \textbf{while} loop for the current value
of $k$. This by itself is not sufficient to ensure that we get to
executing \textbf{set} $\boldsymbol{x}^{k+1}=\mathbf{P}_{T}\boldsymbol{x}^{k,N}$
followed by executing $\mathbf{set}$ $k=k+1$, we must also make
sure that the execution of the code inserted between the lines \textbf{while}
$n<N$ and \textbf{set} $\boldsymbol{x}^{k+1}=\mathbf{P}_{T}\boldsymbol{x}^{k,N}$
terminates in a finite number of steps, for all $k\geq0$.

If the code for the superiorized version of algorithm $\mathbf{P}$
is designed in accordance with the discussion above, then its execution
will produce an infinite sequence $S=\left(\boldsymbol{x}^{k}\right)_{k=0}^{\infty}$
of elements of $\Omega$. However, there is nothing in what we said
so far that indicates anything like ``whenever $O\left(T,\varepsilon,R\right)$
is defined, then $O\left(T,\varepsilon,S\right)$ is also defined
and $\phi\left(O\left(T,\varepsilon,S\right)\right)\leq\phi\left(O\left(T,\varepsilon,R\right)\right)$.''
In order to obtain such a result we need to bring the target function
into the code inserted between the lines \textbf{while} $n<N$ and
\textbf{set} $\boldsymbol{x}^{k+1}=\mathbf{P}_{T}\boldsymbol{x}^{k,N}$.
The way this is typically done (see \citep{CGHH19}) is to break up
the computation from $\boldsymbol{x}^{k,n}$ in $\Delta$ to an $\boldsymbol{x}^{k,n+1}$
in $\Delta$ into stages: first we compute from $\boldsymbol{x}^{k,n}$
a provisional candidate $\boldsymbol{z}$ and check whether $\boldsymbol{z}\in\Delta$
and compare $\phi\left(\boldsymbol{z}\right)$ with previously obtained
$\phi\left(\boldsymbol{x}^{k,n}\right)$. In case of a satisfactory
outcome, we set $\boldsymbol{x}^{k,n+1}$ to $\boldsymbol{z}$, otherwise
we compute a new provisional candidate $\boldsymbol{z}$. The details
of how to do this depend on the definition of $\phi$ and we need
to give up some generality. We do this in the next section.

\section{Perturbation Resilience and Superiorization of Algorithms}

The definitions until now have been very general; in particular, no
restrictions have been placed on the nature of $\Delta$. From now
on we assume, without stating this again, that $\Delta\subseteq\mathbb{R}^{J}$,
for some positive integer $J$ . Thus we are working in a $J$-dimensional
vector space. We use $\left\Vert \boldsymbol{x}\right\Vert $ to denote
the 2-norm of $\boldsymbol{x}\in\Delta$.

An algorithm\emph{ }\textbf{$\mathbf{P}$ }for a problem structure
$\left\langle \Omega,\mathbb{T},\mathcal{P}r\right\rangle $ and for
a set $\Delta\subseteq\mathbb{R}^{J}$ is said to be \textit{strongly
perturbation resilient} \citep{Her12} if, for all $T\in\mathbb{T}$, 
\begin{enumerate}
\item there exists an $\varepsilon\in\mathbb{R}_{+}$ such that $O\left(T,\varepsilon,\left(\left(\mathbf{P}_{T}\right)^{k}\boldsymbol{x}\right)_{k=0}^{\infty}\right)$
is defined for every $\boldsymbol{x}\in\Omega$; 
\item for all $\varepsilon\in\mathbb{R}_{+}$ such that $O\left(T,\varepsilon,\left(\left(\mathbf{P}_{T}\right)^{k}\boldsymbol{x}\right)_{k=0}^{\infty}\right)$
is defined for every $\boldsymbol{x}\in\Omega$, we also have that
$O\left(T,\varepsilon',S\right)$ is defined for every $\varepsilon'>\varepsilon$
and for every sequence $S=\left(\boldsymbol{x}^{k}\right)_{k=0}^{\infty}$
of elements of $\Omega$ with $\boldsymbol{x}^{0}$ arbitrary and
generated by 
\begin{equation}
\boldsymbol{x}^{k+1}=\mathbf{P}_{T}\left(\boldsymbol{x}^{k}+\beta_{k}\boldsymbol{v}^{k}\right),\:\mathrm{for\: all\:\,}k\geq0,\label{eq:perturbations}
\end{equation}
where $\beta_{k}\boldsymbol{v}^{k}$ are \textit{bounded perturbations},
meaning that the sequence $\left(\beta_{k}\right)_{k=0}^{\infty}$
of nonnegative real numbers is \textit{summable} (that is, ${\displaystyle \sum\limits _{k=0}^{\infty}}\beta_{k}\,<\infty$),
the sequence $\left(\boldsymbol{v}^{k}\right)_{k=0}^{\infty}$ of
vectors in $\mathbb{R}^{J}$ is bounded and, for all $k\geq0$, $\boldsymbol{x}^{k}+\beta_{k}\boldsymbol{v}^{k}\in\Delta$. 
\end{enumerate}
From the practical point of view, it is important to note that, even
though the definition of strong perturbation resilience appears to
be complicated (and hence hard to check), easily verifiable sufficient
conditions for strong perturbation resilience have been published
\citep[Theorem 1]{Her12}.

We now discuss how this specific definition relates to the general
discussion of the previous section.

Note that results proved for strongly perturbation resilient algorithms
are by necessity weaker than the ``rather-strong mathematical restatement''
prior to the pseudocode in the previous section. One reason for this
is that not all algorithms are strongly perturbation resilient. For
example, the first condition above for an algorithm to be strongly
perturbations resilient is not automatically satisfied by all algorithms,
even though one may argue that most algorithms of practical interest
will satisfy that condition.

Definitions of the sequence $S$ in the definition above and prior
to the pseudocode in the previous section can be made to match up
by selecting $\mathbf{S}$ so that $\mathbf{S}\boldsymbol{x}^{k}=\boldsymbol{x}^{k}+\beta_{k}\boldsymbol{v}^{k}$;
see (\ref{eq:perturbations}). Nevertheless, even for strongly perturbation
resilient algorithms, we cannot conclude that ``whenever $O\left(T,\varepsilon,R\right)$
is defined, then $O\left(T,\varepsilon,S\right)$ is also defined''
because the definition above guarantees only something weaker, namely
that whenever ``$O\left(T,\varepsilon,\left(\left(\mathbf{P}_{T}\right)^{k}\boldsymbol{x}\right)_{k=0}^{\infty}\right)$
is defined for every $\boldsymbol{x}\in\Omega$, we also have that
$O\left(T,\varepsilon',S\right)$ is defined for every $\varepsilon'>\varepsilon$.''
However, from the practical point of view, we are very near to where
we wanted to be, since the condition ``$O\left(T,\varepsilon,\left(\left(\mathbf{P}_{T}\right)^{k}\boldsymbol{x}\right)_{k=0}^{\infty}\right)$
is defined for every $\boldsymbol{x}\in\Omega$'' would be satisfied
by most commonly-used algorithms that we may wish to superiorize and
``$O\left(T,\varepsilon',S\right)$ is defined for every $\varepsilon'>\varepsilon$''
is in practice as useful as ``$O\left(T,\varepsilon,S\right)$ is
defined.''

Next we discuss how to design the computation from $\boldsymbol{x}^{k,n}$
in $\Delta$ to an $\boldsymbol{x}^{k,n+1}$ in $\Delta$ (see the
last paragraph of the previous section) so that the resulting $\boldsymbol{x}^{k,N}$
will be of the form $\boldsymbol{x}^{k}+\beta_{k}\boldsymbol{v}^{k}$
with bounded perturbations $\beta_{k}\boldsymbol{v}^{k}$, see the
lines following (\ref{eq:perturbations}). This is done by specifying
a summable sequence $\left(\gamma_{\ell}\right)_{\ell=0}^{\infty}$
of nonnegative real numbers (for example, $\gamma_{\ell}=a^{\ell}$,
where $0<a<1$). We set the value of $\ell$ to -1 in the line(s)
to be inserted above \textbf{repeat} in the skeleton of the superiorized
algorithm and make sure that $\ell$ is increased by 1 prior to each
one of its subsequent uses in the superiorized algorithm. (This and
what is described in the next sentence are done, for example, in the
two superiorized versions, called Algorithms 1 and 2, in \citep{CGHH19}.)
The provisional candidate $\boldsymbol{z}$ is defined as $\boldsymbol{z}=\boldsymbol{x}^{k,n}+\gamma_{\ell}\boldsymbol{v}^{k,n}$,
where $\left\Vert \boldsymbol{v}^{k,n}\right\Vert \leq1$. Also, the
$\boldsymbol{v}^{k,n}$ need to be selected so that $\boldsymbol{x}^{k,n+1}$
is obtained from $\boldsymbol{x}^{k,n}$ in a finite number of steps;
implying that only finitely many provisional candidates need to be
tried before finding one that is acceptable as the choice of $\boldsymbol{x}^{k,n+1}$;
examples for this are provided in the literature \citep{Cen19}, in
particular in \citep{Her12,CGHH19}. That the resulting $\boldsymbol{x}^{k,N}$
will be of the form $\boldsymbol{x}^{k}+\beta_{k}\boldsymbol{v}^{k}$
with bounded perturbations $\beta_{k}\boldsymbol{v}^{k}$ follows
from a discussion on p. 5538 of \citep{Her12}.

The code for calculating the $\boldsymbol{v}^{k,n}$ also needs to
be inserted into the skeleton between the lines \textbf{while} $n<N$
and \textbf{set} $\boldsymbol{x}^{k+1}=\mathbf{P}_{T}\boldsymbol{x}^{k,N}$.
The purpose of this code should be to find a provisional candidate
$\boldsymbol{z}$ with a reduced target value. Such code depends on
the definition of $\phi$. For example, if we can calculate some partial
derivatives of $\phi$ at $\boldsymbol{x}^{k,n}$, then they can be
used to produce suitable provisional candidates. Alternatively, one
can just use a component-wise search in the neighborhood of $\boldsymbol{x}^{k,n}$
for finding provisional candidates. Details of these two approaches
are given in Algorithms 1 and 2, respectively, of \citep{CGHH19}.

\section{Concluding Discussion}

The superiorization methodology has been found efficacious for solving
real world problems in many areas of application \citep{Cen19}. This
brief paper attempts to get to that bare essence of the methodology
which is the common core in past, present and future applications.
For this an abstract approach has been adopted, in the spirit of the
study of the general concept of algebraic structures. Specifically,
targeted problem structures $\left\langle \Omega,\mathbb{T},\mathcal{P}r,\phi\right\rangle $
have been introduced, where:
\begin{itemize}
\item $\Omega$ is an arbitrary nonempty set. It is the domain of the problems
that we wish to ``solve.''
\item $\mathbb{T}$ is another arbitrary nonempty set. Its element are used
to index the problems that we wish to ``solve.'' Hence $\mathbb{T}$
is referred to as a problem set.
\item $\mathcal{P}r$ is a called the proximity function, it is used to
specify the problems $T$ in $\mathbb{T}$. For any $T\in\mathbb{T}$,
$\mathcal{P}r_{T}$ maps $\Omega$ into \textcolor{black}{$\mathbb{R}_{+}$.
The intuitive idea is that for any $\boldsymbol{x}\in\Omega$,$\;\,$
}$\mathcal{P}r_{T}(\boldsymbol{x})$ indicates the undesirability
of $\boldsymbol{x}$ as a solution of the problem $T$.
\item $\phi$\textcolor{black}{{} assigns a real value to every element x
of \textgreek{W};$\;\,$ $\phi\left(\boldsymbol{x}\right)$ indicates
the prior undesirability of $\boldsymbol{x}$ without considering
any specific problem }$T\in\mathbb{T}$.
\end{itemize}
The aim of the superiorization methodology is to produce for an algorithm
a version of it whose performance is as good as that of the original
algorithm from the point of view of the proximity function, but is
better from the point of view of the target function. We can make
this more precise as in the following description.

Let\textbf{ $\mathbf{P}$ }be an algorithm for a problem structure
$\left\langle \Omega,\mathbb{T},\mathcal{P}r\right\rangle $ and a
given $\Delta$. For a $T\in\mathbb{T}$ and $\varepsilon\in\mathbb{R}_{+}$,
assume that there is a first element $O\left(T,\varepsilon,R\right)$
in the sequence $R=\left(\boldsymbol{x}^{k}\right)_{k=0}^{\infty}$
produced by repeated application of $\mathbf{P}_{T}$ for which $\mathcal{P}r_{T}(\boldsymbol{x})\leq\varepsilon$.
A desirable aim of superiorization is that if the superiorized version
of algorithm $\mathbf{P}$ is initialized with $\boldsymbol{x}^{0}$,
then it produces a sequence $S$ of elements of $\Omega$ such that
$\phi\left(O\left(T,\varepsilon,S\right)\right)\leq\phi\left(O\left(T,\varepsilon,R\right)\right)$.
This has not yet been achieved, but much progress has been made toward
it.

\end{document}